\documentclass[twoside,british]{article}
\pdfoutput=1
\usepackage[T1]{fontenc}
\usepackage[latin9]{inputenc}
\usepackage{amsmath}
\usepackage{amssymb}
\usepackage{graphicx}

\makeatletter
\newcommand{\lyxaddress}[1]{
\par {\raggedright #1
\vspace{1.4em}
\noindent\par}
}

\@ifundefined{date}{}{\date{}}
\pdfoutput=1

\@ifundefined{showcaptionsetup}{}{%
 \PassOptionsToPackage{caption=false}{subfig}}
\usepackage{subfig}
\makeatother

\usepackage{babel}
\begin{document}

\title{Levy process simulation by stochastic step functions}

\author{Torquil Macdonald S{\o}rensen%
\thanks{Email: t.m.sorensen@matnat.uio.no, torquil@gmail.com%
} \and Fred Espen Benth%
\thanks{Email: fredb@math.uio.no%
}}

\maketitle

\lyxaddress{\begin{center}
Centre of Mathematics for Applications\\
University of Oslo, NO-0316 Oslo, Norway
\par\end{center}}
\begin{abstract}
We study a Monte Carlo algorithm for simulation of probability distributions
based on stochastic step functions, and compare to the traditional
Metropolis/Hastings method. Unlike the latter, the step function algorithm
can produce an uncorrelated Markov chain. We apply this method to
the simulation of Levy processes, for which simulation of uncorrelated
jumps are essential.

We perform numerical tests consisting of simulation from probability
distributions, as well as simulation of Levy process paths. The Levy
processes include a jump-diffusion with a Gaussian Levy measure, as
well as jump-diffusion approximations of the infinite activity NIG
and CGMY processes.

To increase efficiency of the step function method, and to decrease
correlations in the Metropolis/Hastings method, we introduce adaptive
hybrid algorithms which employ uncorrelated draws from an adaptive
discrete distribution defined on a space of subdivisions of the Levy
measure space.

The nonzero correlations in Metropolis/Hastings simulations result
in heavy tails for the Levy process distribution at any fixed time.
This problem is eliminated in the step function approach. In each
case of the Gaussian, NIG and CGMY processes, we compare the distribution
at $t=1$ with exact results and note the superiority of the step
function approach.
\end{abstract}

\section{Introduction}

Levy processes are a type of stochastic process whose paths can behave
erratically like a Brownian motion, as well as include discontinuities,
i.e. jumps. Many observed phenomena in nature and human society display
continuous erratic motion similar to Brownian motion, while also having
random jumps at random times, and can therefore be modelled using
Levy processes. Areas of application include e.g. microscopic physics,
chemistry, biology and financial markets. Therefore the study of computer
simulation methods for Levy processes is an important subject.

In our case, will will employ a jump-diffusion approximation when
we apply our method to infinite activity Levy processes. Jump-diffusion
processes can be considered as a sum of three independent component
processes. A deterministic drift, a Brownian motion, and lastly a
finite activity jump process. The jump process is completely described
by a Levy measure $\nu$ on $\mathbb{R}$. The average jump rate,
or intensity, is given by the total weight $\lambda:=\nu(\mathbb{R})$,
and the distribution of jump values follows the Levy probability measure
$\nu_{1}:=\nu/\lambda$.

By definition, the Brownian motion has independent increments, and
this is also the case for subsequent jumps in the jump process. In
computer simulations, violation of these properties will give incorrect
results. Simulation of the Brownian motion is easy, since well known
algorithms exist for producing uncorrelated draws from a normal distribution,
using a good pseudo-random number generator (PRNG).

On the other hand, the simulation of the jump process requires more
care. In principle we can obtain uncorrelated draws from $\nu_{1}$
by inverting its cumulative distribution function. However, this might
not be feasible to do for a given $\nu_{1}$. In this case, the Metropolis/Hastings
(MH) algorithm might come to the rescue. It can easily produce a Markov
chain with values distributed according to $\nu_{1}$. However, subsequent
values will be correlated, as described below.

In this paper we propose an algorithm to produce uncorrelated jumps
along each path, without generating such a multitude of Markov chains.
This method is based on stochastic step functions (SF), which will
be defined below. As opposed to MH, it is not an accept/reject algorithm.
Therefore it is able to generate an uncorrelated chain of values distributed
according to any given probability measure. We will test this algorithm
in jump-diffusion computer simulations and compare with MH. We apply
the simulation techniques to a Gaussian process, for which the exact
distribution at $t=1$ is known. Convergence towards this exact result
is studied.

We also consider adaptive variants of these algorithms. For these,
we subdivide $\mathbb{R}$ into appropriate regions, and generate
a discrete distribution that allows us to efficiently draw among these
regions in an uncorrelated way. As the simulation progresses, this
discrete distribution is adaptively improved. The calculation of the
jump rate $\lambda$ is also adaptively improved.

As an interesting application of these simulation techniques, we also
look at infinite activity pure jump processes, which are also a Levy
process subclass. Here, motion occurs in the form of an infinitude
of discontinuous jumps. Some of these processes can be approximated
by jump-diffusion by substituting the smallest jumps for an appropriate
Brownian motion \cite{AsmussenRosinski:2001}. The examples we focus
on are the NIG and CGMY processes.

In section \ref{sec:JumpDiff} we review the elementary facts about
jump-diffusion processes. In section \ref{sec:ProbMeasSim} we describe
the relevant simulation methods, and discuss their pros and cons,
as well as provide results from numerical experiments. In section
\ref{sec:SubdivAdapt} we describe our probability space subdivision
method and the accompanying adaptability properties of the algorithms,
and study the efficiency and correlation strengths of different algorithms
by simulation from a Gaussian distribution. We then proceed in section
\ref{sec:CompSimJumpDiff} to perform simulations of jump-diffusion
processes. We notice how the MH correlations adversely affect the
distribution of the simulated process, and compare the SF and MH algorithms
for simulation of a process with a Gaussian Levy measure. Lastly,
in section \ref{sec:SimInfAct}, we review the jump-diffusion approximation
of infinite activity pure jump processes, and apply our simulations
techniques on the infinite activity NIG and CGMY processes, in order
to compare algorithms in these interesting cases.

\section{Jump-diffusion processes\label{sec:JumpDiff}}

First we review some elementary facts about real-valued jump-diffusion
Levy processes on a time interval $[0,T]$. Such a process can be
expressed as a sum of three simple components

\begin{equation}
L_{t}=\mu t+B_{t}+J_{t}.\label{eq:JumpDiffDecomp}
\end{equation}
The first part is a deterministic drift with rate $\mu$, $B_{t}$
is Brownian motion, and $J_{t}$ is a compound Poisson process. The
latter describes completely the discontinuities (jumps) of the paths
of $L_{t}$, by means of a \emph{Levy measure} $\nu$ on $\mathbb{R}$.
We will several times abuse notation by writing $\nu$ both for the
Levy measure and its density. Firstly, this measure determines the
jump intensity (jump rate)
\begin{equation}
\lambda:=\nu(\mathbb{R})<\infty,\label{eq:intensity}
\end{equation}
 of $L_{t}$. Secondly, the corresponding \emph{Levy probability measure}

\begin{equation}
\nu_{1}:=\nu/\lambda,\label{eq:LevyProbMeas}
\end{equation}
determines the distribution of jump sizes. Note that we do not have
$E[L_{t}]=\mu t$ in general, although this is satisfied in our examples.

Now, $J_{t}$ can be expressed as

\begin{equation}
J_{t}=\sum_{j=1}^{N_{t}}V_{j},\label{eq:JumpProc}
\end{equation}
where $N_{t}$ is a Poisson process of intensity $\lambda$, and the
jumps $\{V_{j}\}$ are indendent random variables distributed according
to the Levy probability measure. We do not simulate the Poisson process
$N_{t}$ directly. Instead, for each path we draw the total number
of jumps individually from an exponential distribution with average
$\lambda T$. Then we randomly distribute those jumps on $[0,T]$.
We will choose $T=1$ for our simulations. The lengths of consecutive
time intervals between jumps will be independent and exponentially
distributed, and give the correct jump intensity. This is described
in \cite{ContTankov:2003}.

Note that for a general Levy process $\lambda$ is in general not
finite, because the Levy measure might diverge too strongly towards
the origin. However, the following condition always holds,

\begin{equation}
\int_{\mathbb{R}}\min(1,s^{2})\nu(ds)<\infty.\label{eq:LevyMeasCond}
\end{equation}
which restricts the strength on the divergence of $\nu$ at the origin.

The difficulty in simulation is to draw independent jump values from
$\nu_{1}$. We will generate a Markov chain with $\nu_{1}$ as its
invariant measure. However, it is essential for correct simulation
that jumps along each path are uncorrelated. Note that existence of
correlations between jumps in \emph{different} paths will not ruin
the convergence in distribution, but only slow it down.

For a Markov chain $\{X_{i}\}$, we define the sequential correlation
as
\begin{equation}
c:=\mathrm{E}[(X_{i+1}-\bar{X})(X_{i}-\bar{X})]/\sigma^{2},\label{eq:correlation}
\end{equation}
where $\bar{X}$ is the Markov chain average and $\sigma^{2}$ is
its standard deviation,
\[
\bar{X}:=E[X],\quad\sigma^{2}:=E[X^{2}-\bar{X}^{2}].
\]

\section{Simulation of a probability measure\label{sec:ProbMeasSim}}

To simulate the jumps, we must draw independent values from $\nu_{1}$.
In cases where $\nu_{1}$ is complicated, it is common to use the
Metropolis/Hastings (MH) algorithm \cite{Metropolis:1953,Hastings:1970}.
This method generates a Markov chain with $\nu_{1}$ as its invariant
density. Unfortunately, the Markov chain often has large correlations
between successive values. Successive values in such a chain cannot
be used to represent jumps within a single jump-diffusion path.

It is possible to reduce this problem by several methods. One is to
skip a number of terms in the Markov chain to reduce correlation.
This results in a loss of efficiency. Another method is to generate
multiple independent Markov chains. Each jump-diffusion path then
uses values from different Markov chains. This will lead to a correct
pathwise simulation, and therefore correct convergence in distribution.
The correlations will in this case only slow down the convergence.

We will look at two methods, MH and one based on stochastic step fuctions
(SF). Both rely on a transition probability distribution $\rho$ to
determine upcoming values from the previous one. We say that we are
dealing with a \emph{local} algorithm if $\rho$ is localized around
the current value. An \emph{independent} algorithm results if $\rho$
is independent of the current position.

Use of an independent MH algorithm can reduce correlations dramatically
if one is able to find a $\rho$ that approximates $\nu_{1}$. The
downside is that efficiency will suffer if $\rho$ is not computationally
simple. We will focus on generation of low correlation Markov chains
in order to get by using only one chain for the Levy process jumps.

\subsection{Metropolis/Hastings (MH)\label{sub:LMH}}

The well-known MH algorithm generates a realisation $\{x_{i}\}$ of
a Markov chain distributed according to $\nu_{1}$. Its most popular
incarnation is local, where the transition probability density $\rho(x_{i},x_{i+1})$
involved in each transition $x_{i}\mapsto x_{i+1}$ is localized around
the value $x_{i}$, and its width is adjusted to give an average acceptance
rate somewhere around $1/2$. The resulting correlation between successive
values can be classified into two types:
\begin{itemize}
\item \emph{Rejection correlation}: Since both the local and independent
algorithms are based on acceptance/rejectance, correlations are introduced
by repeated values due to rejections. With an acceptance rate around
$1/2$, repeated values will often occur.
\item \emph{Transition correlation}: In order for the local algorithm to
be efficient, its transition probability measure $\rho(x_{i},x_{i+1})$
must in many cases be quite strongly localized; otherwise too small
an acceptance rate will result. Thus the difference between subsequent
variates of the Markov chain will tend to be small.
\end{itemize}
To reduce transition correlation, the width of the transition distribution
can be increased. But this typically leads to a reduced acceptance
ratio, which gives an increased rejection correlation, and vice versa.

\begin{figure}
\centering{}\includegraphics[width=8cm]{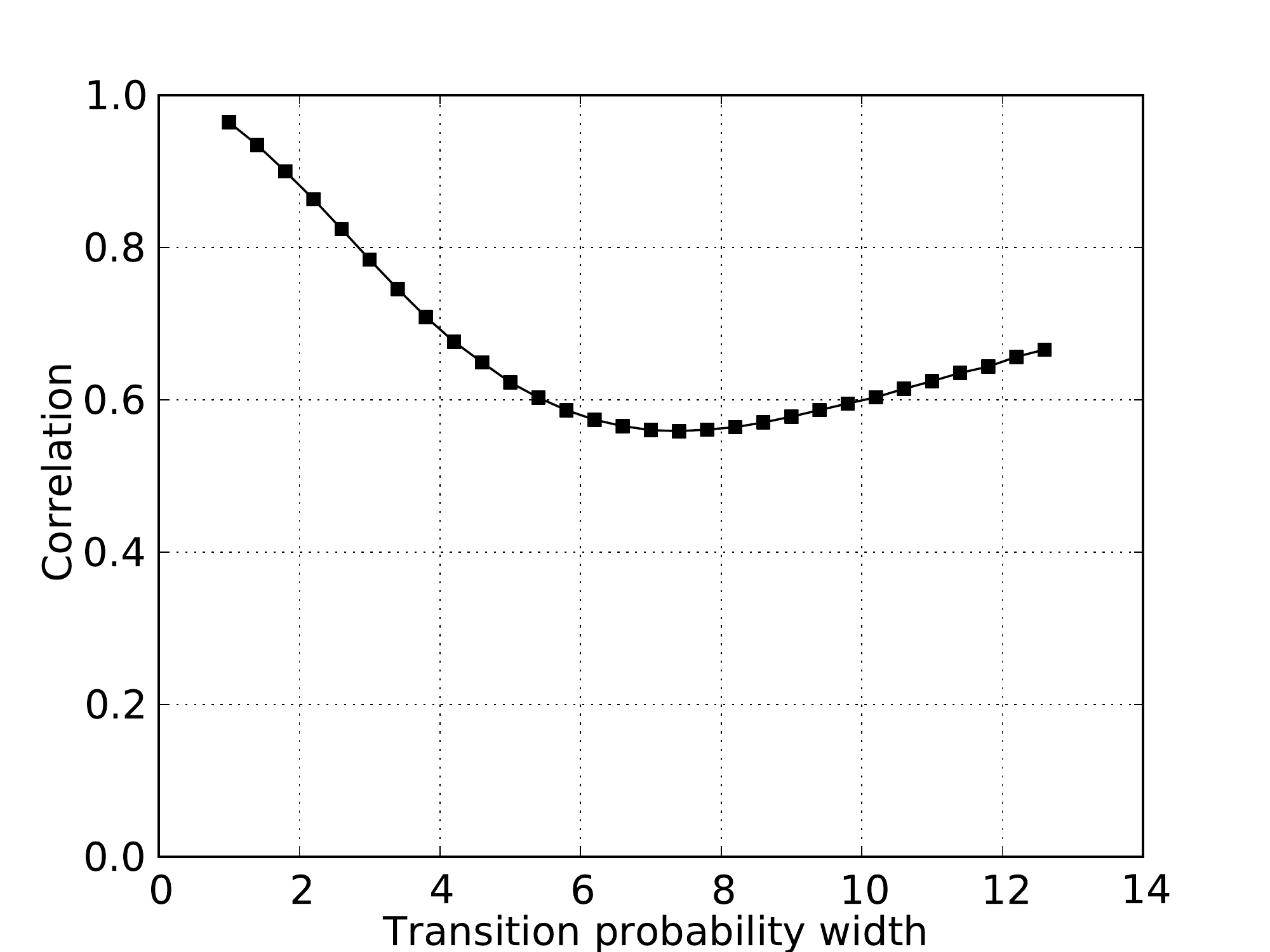}\caption{Correlation $c$ versus transition probability width $w$ for the
local MH algorithm, in a simulation of the unit normal distribution
using a uniform centered transition probability distribution. Correlation
is minimal around $w\approx7$, where $c\approx0.55$.\label{fig:CorrVsWidthLMH}}
\end{figure}

As a numerical example, consider a unit variance normal distribution,
simulated with a simple localized uniform transition probability measure.
The correlation defined in \eqref{eq:correlation} as a function of
transition probability width is displayed in Figure \ref{fig:CorrVsWidthLMH}.
As expected, the correlation has a minimum. Towards smaller widths
the transition correlation increases, and towards larger widths the
rejection correlation increases. Thus there is a lower bound on the
amount of correlation for this algorithm, and it seems that the local
MH algorithm is therefore not suited for our purposes. We will from
now on focus on the independent MH algorithm.

\subsection{Stochastic step function (SF)\label{sub:SF}}

Next, we propose a simulation algorithm based on stochastic step functions,
that can be adjusted to completely eliminate correlations. This possibility
of vanishing correlations is an attractive property that allows it
to be used as a reference algorithm. It can also be adjusted to run
more efficiently, with nonzero correlation.

Let us give a pathwise definition of the SF process on $[0,\infty)$
for some measure density $\nu$ on $\Omega\subset\mathbb{R}$, not
necessarily normalized. Consider a sequence $\{\tau_{i}\}_{i=0}^{\infty}\subset(0,\infty)$
of \emph{resting times} and corresponding \emph{jump times} $\{t_{i}\}_{i=0}^{\infty}$,
defined recursively by
\[
t_{0}=0,\quad t_{i+1}:=t_{i}+\tau_{i}.
\]
In addition, let $\{s_{i}\}_{i=0}^{\infty}\subset\Omega$ be a Markov
chain with initial value $s_{0}=0$ and transition probability density
$\rho(s_{i},s_{i+1})$. Assume that $\rho(s_{i},\cdot)$ is absolutely
continuous with respect to the Lebesgue measure on $\mathbb{R}$,
and homogeneous, i.e. it can be expressed as $\rho(s_{i+1}-s_{i})$.

From this data, we are now ready to define our stochastic step function
process $X_{t}$, by defining its paths as the piecewise constant
\emph{cadlag} functions (right-continuous with left limits) of the
form
\[
X_{t}=\sum_{i=1}^{\infty}s_{i}\chi_{I_{i}}(t),\quad I_{i}:=[t_{i},t_{i+1}),
\]
where $\chi_{I}$ is the indicator function for the interval $I$,
and where the resting times are given by $\tau_{i}:=\nu(s_{i})$.

We see that the paths of $X_{t}$ linger for a some time at each of
its attained positions, with a resting time equal to the local value
of the density $\nu$. Consider the graph of such a path on $[0,T]$
for large $T$ compared to $\sup\nu$. As $T$ increases, the relative
path length within a given subset $A\subset\Omega$ on the vertical
axis converges towards $\nu(A)/\nu(\Omega)=:\nu_{1}(A)$. When sampling
this path on a uniform time-grid, the obtained values will therefore
be distributed according to the probability density $\nu_{1}$.

As in the case of local MH, we get transition probability correlations
if $\rho$ is localized. In this algorithm, however, there is no accept/reject
step, and therefore no rejection correlation. As an example, let us
choose $\rho$ to be the uniform probability distribution on $\Omega$.
As noted above, if we sample the step function path generated above
on a uniform time grid, we get a chain of values distributed according
to $\nu_{1}$. If we make sure that the time discretization interval
size $\Delta t$ is larger than $\sup\nu$, repeated values will never
occur and there are no correlations.

Note that as in MH algorithms, we do not need to know scale factor
relating $\nu_{1}$ and $\nu$. If $\sup\nu$ is initially unknown,
we can simply start with any estimate, and improve it adaptively as
the algorithm runs.

It is easy to see that the amount of computer time spent between each
discrete time value is proportional to $\sup\nu/\nu(\Omega)$ in its
correlation-free mode. If this ratio is large, the algorithm will
be inefficient.

One can make the time discretization finer to increase the rate of
variate generation, but this introduces correlations since values
for which $\nu$ is large will be repeated with complete certainty.
This is similar to the Metropolis algorithm, where values of large
$\nu$ are repeated (although not with complete certainty). The SF
algorithm can be made less deterministic in this case by letting the
resting times be random variables distributed according to an exponential
distribution with mean $\nu(x)$, where $x$ is the current position.

Since we are concerned with minimizing correlations in the context
of jump-diffusion simulations, we will use independent algorithms,
where $\rho(x_{i},x_{i+1})$ is independent of $x_{i}$.

\subsection{Simulation comparison}

Let us illustrate the advantages of the local SF algorithm over local
MH by considering an example with a probability distribution with
several modes on a sample space $\Omega=[0,1]$. Our goal is to simulate
values from a probability distribution proportional to the following
unnormalized density with two strong modes,
\begin{equation}
\nu(x)=\left\{ \begin{array}{ll}
1 & ,x\in[0,0.25)\cup[0.5,0.75)\\
0.01 & ,\mbox{otherwise}.
\end{array}\right.\label{eq:TwoModeProbDens}
\end{equation}
For both MH and SF, we used a localized uniform transition probability
density $\rho$ of width $1/2$. For MH, this gave an acceptance rate
of approximatively 0.55. The SF Markov chain realisation was obtained
by sampling the stochastic step function using a lattice spacing slightly
larger than $\sup\nu$ in order to avoid any repeating values.

\begin{figure}
\begin{centering}
\subfloat[Local MH]{\includegraphics[width=6cm]{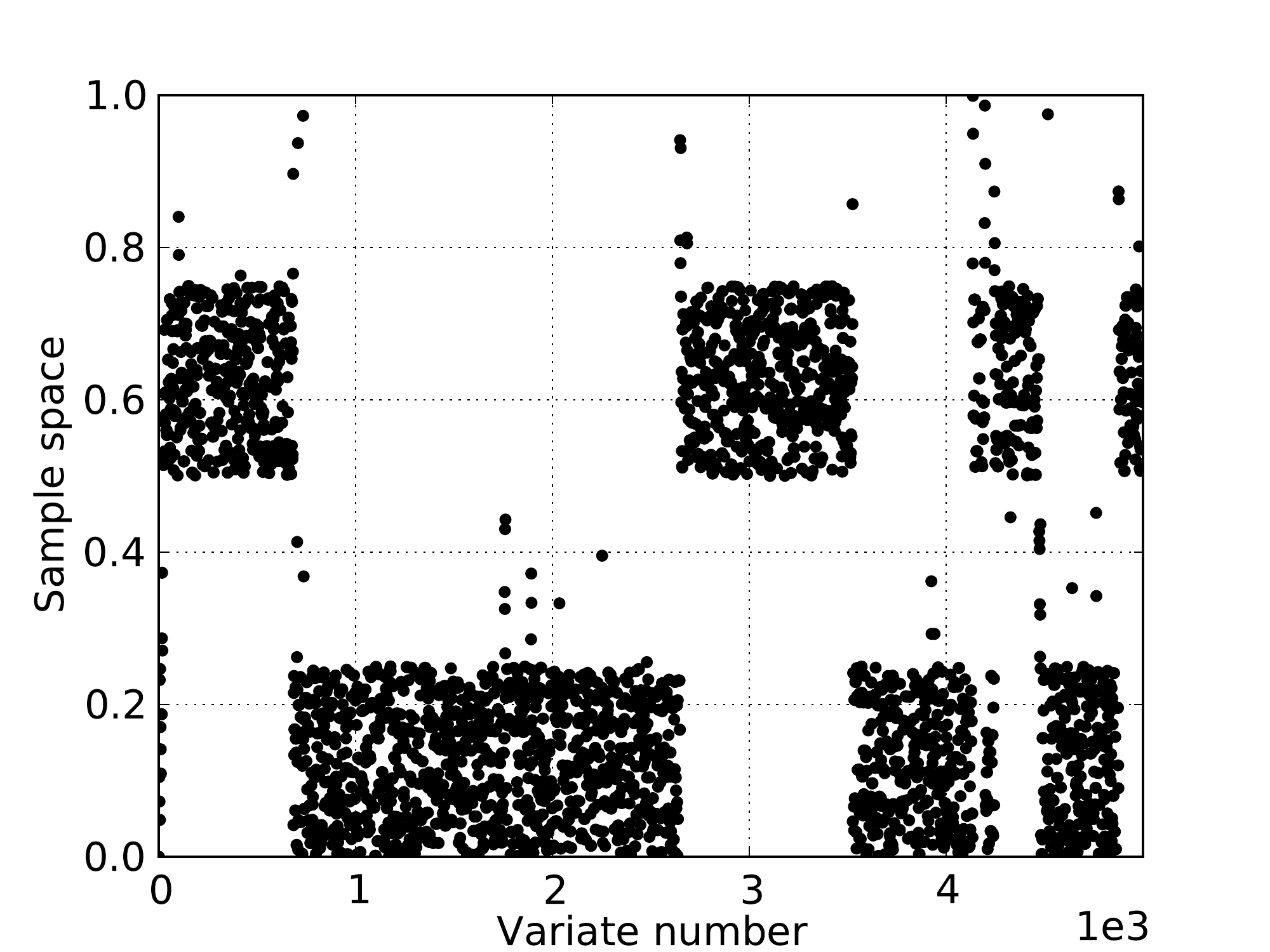}

}\subfloat[Local SF]{\includegraphics[width=6cm]{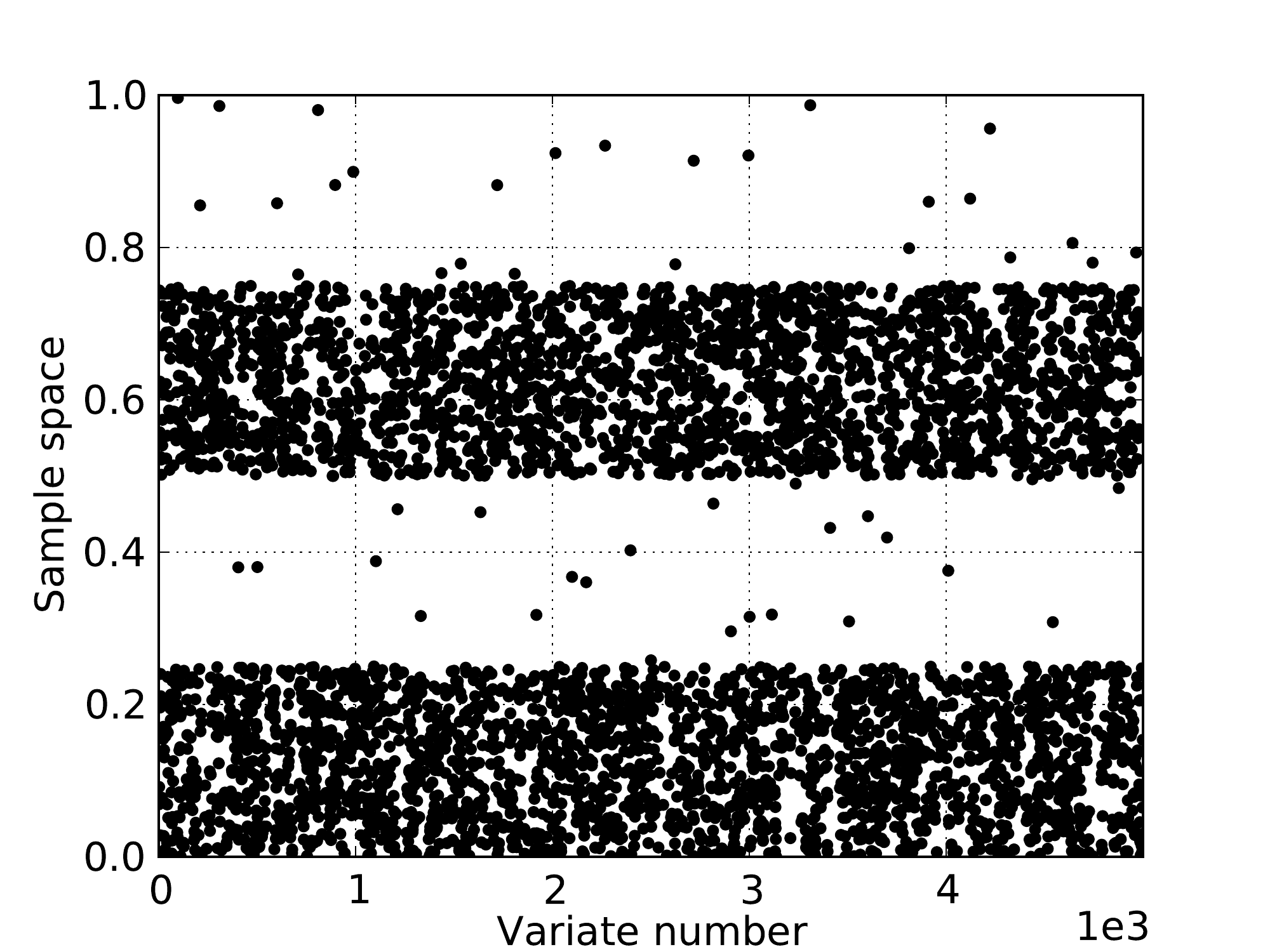}

}\caption{Markov chain realisations using local MH and local SF. The SF algorithm
mixes much better than MH.\label{fig:Modes_MH_vs_SF}}

\par\end{centering}

\end{figure}

The results are shown in Figure \ref{fig:Modes_MH_vs_SF}. One sees
immediately that the MH algorithm has the potential of getting stuck
inside a mode for long periods. This is caused by rejection correlation
as found in accept/reject algorithms such as MH. On the other hand,
the SF algorithm has no such problems since it lacks such correlations.
This serves to illustrate a problem that often can occur with MH.
The SF algorithm is a clear favourite in this case if mixing is important
for the application of these Markov chains.

\section{Probability space subdivision and adaptability\label{sec:SubdivAdapt}}

The Levy probability measure from which we will simulate in the context
of jump-diffusions will often be defined on the unbounded probability
space $\Omega=\mathbb{R}$. This presents no difficulty for the MH
algorithm. SF algorithms on the other hand need to impose upper and
lower cutoffs on $\Omega$, in order to avoid step functions that
diverge towards infinity. For simplicity, we impose such cutoffs on
both algorithms in our examples. Alternatively, one could perform
a topology-changing coordinate transformation on $\Omega$ to obtain
a compact space, however we will not do this in our examples. From
now on we therefore assume $\Omega\subset\mathbb{R}$ to be a bounded
interval which we choose symmetrically about the origin. We will only
deal with symmetrical Levy measures in our examples.

In order to reduce correlations in the MH algorithms, and increase
efficiency in the SF algorithms, we define a finite disjoint subdivision
$\{U_{i}\}$ of $\Omega$, where $\cup_{i}U_{i}=\Omega$. We construct
a discrete distribution $\tilde{\nu}$ on the finite set $\{U_{i}\}$.
For each $U_{i}\subset\Omega$, the value $\tilde{\nu}(U_{i})$ must
approximate $\nu(U_{i})$, i.e. the Levy measure weight of $U_{i}$.
The simulation algorithm starts with a preliminary estimate of $\tilde{\nu}$
which is adaptively improved throughout the simulation, as explained
below in the two cases of MH and SF. This is reminiscent of variance
reduction schemes used in Monte Carlo integration, such as stratification
and the VEGAS algorithm \cite{Lepage1978}. We now describe in more
details how this is done in each case.

\subsection{Adaptive Independent MH (AIMH)\label{sub:AIMH}}

As explained above, we will use an independent MH algorithm. The independent
transition probability $\rho$ is defined as follows. First, use the
discrete unnormalized probability measure $\tilde{\nu}$ to draw a
subset $U_{i}$. Then a random position within this subdomain is proposed
and the value of $\nu$ at this position is calculated. Thereafter
follows the usual MH accept/reject step.

The initial draw of $U_{i}$ is done without introducing correlations,
using e.g. an efficient algorithm which is independent of the normalization
of $\tilde{\nu}$ \cite{Walker:1977}. The registered values of $\nu$
are accumulated, and used periodically in the simulation to improve
$\tilde{\nu}$. Essentially, a separate Monte Carlo simulation is
being performed within each subdomain $U_{i}$ to improve the discrete
distribution $\tilde{\nu}$, while the algorithm proceedes to generate
new draws from $\nu$.

Since subdomains $U$ of $\Omega$ are drawn by means of $\tilde{\nu}$
without introducing correlations, the amount of correlation generated
by the algorithm as a whole is reduced. Transition correlation is
completely eliminated since since the algorithm is independent. It
is impossible to completely eliminate rejection correlation in an
MH algorithm unless $\rho$ is identical to $\nu_{1}$. However, since
$\nu$ is well approximated on each subdomain (as long as the subdivision
is sufficiently fine), these will be small, and the acceptance rate
will be high. As $\tilde{\nu}$ adapts, the acceptance rate rises
and correlation decreases. The subdivision discretization implies
a lower bound for the amount of correlation. But this is still a big
improvement on a basic local MH algorithm for the cases we consider.
Note that the subdivision cannot be made arbitrarily fine, since the
discrete algorithm will then become inefficient.

\subsection{Adaptive Independent SF (AISF)\label{sub:AISF}}

It is possible to improve the efficiency of the SF algorithm by a
similar construction, without introducing correlations. In this case,
the SF algorithm proceeds as follows. First draw a subdomain $U_{i}$
using $\tilde{\nu}$. Draw a random position $x$ within this subdomain,
and record the value $\nu(x)$ for later use to improve $\tilde{\nu}$.
For each subset $U_{i}$ we keep an estimate of $\sup_{i}\nu$ which
is updated for each sample. Each subset is also accompanied by a local
time variable $t_{i}$. For each position $x$ within $U_{i}$, we
increase $t_{i}$ by the resting time $\nu(x)$. We choose new positions
independently within $U_{i}$ until $t_{i}$ exceeds the current estimate
of $\sup_{i}\nu$.

When the above loop exists, we have determined our new position and
we subtract $\sup_{i}\nu$ from $t_{i}$ in anticipation of the next
time we enter $U_{i}$. In a sense we are using a different time discretization
in each subset $U_{i}$ of $\Omega$.

This modified algorithm avoids spending a lot of time in areas of
low probability for two reasons. First, the low probability subdomains
will rarely be entered when drawing from the dicrete distribution
$\tilde{\nu}$. Secondly, the amount of time necessary to exit the
loop in those regions is much smaller, since the local $\sup\nu$
is small. By exploiting the subdivision of $\Omega$, and using an
efficient algorithm for drawing from discrete distributions, we improve
the efficiency of the SF algorithm drastically in nontrivial cases.

Both $\tilde{\nu}$ and the estimates $\sup\nu$ are adaptively improved
throughout the simulation.

\subsection{Simulation test}

To check the rate of convergence of the differerent methods, we simulated
from a Gaussian distribution with unit variance. The rate of variate
generation, correlation, and the convergence towards the known exact
result were analyzed.

\begin{figure}
\begin{centering}
\subfloat[Local MH]{\includegraphics[width=6cm]{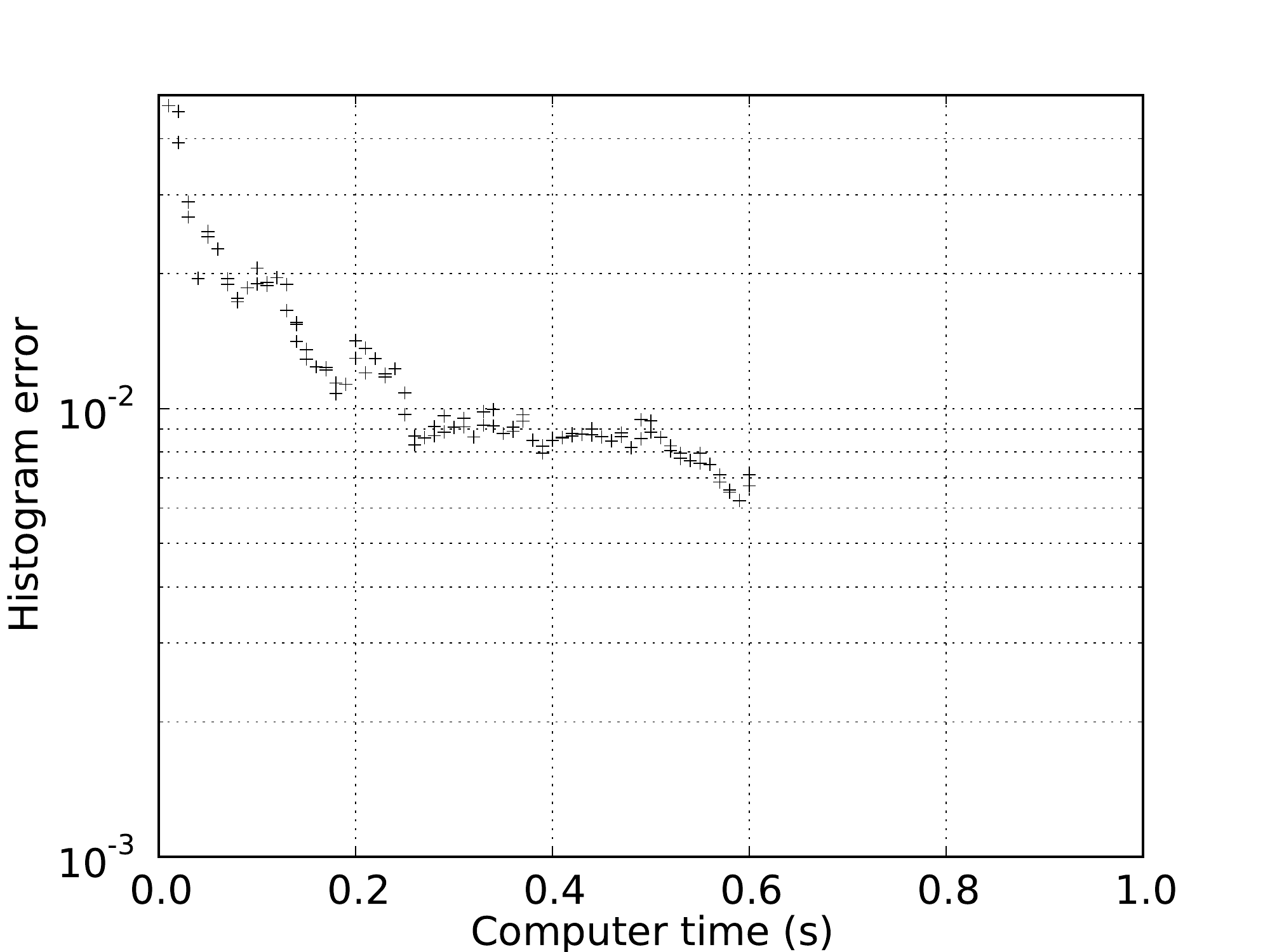}

}\subfloat[AIMH]{\includegraphics[width=6cm]{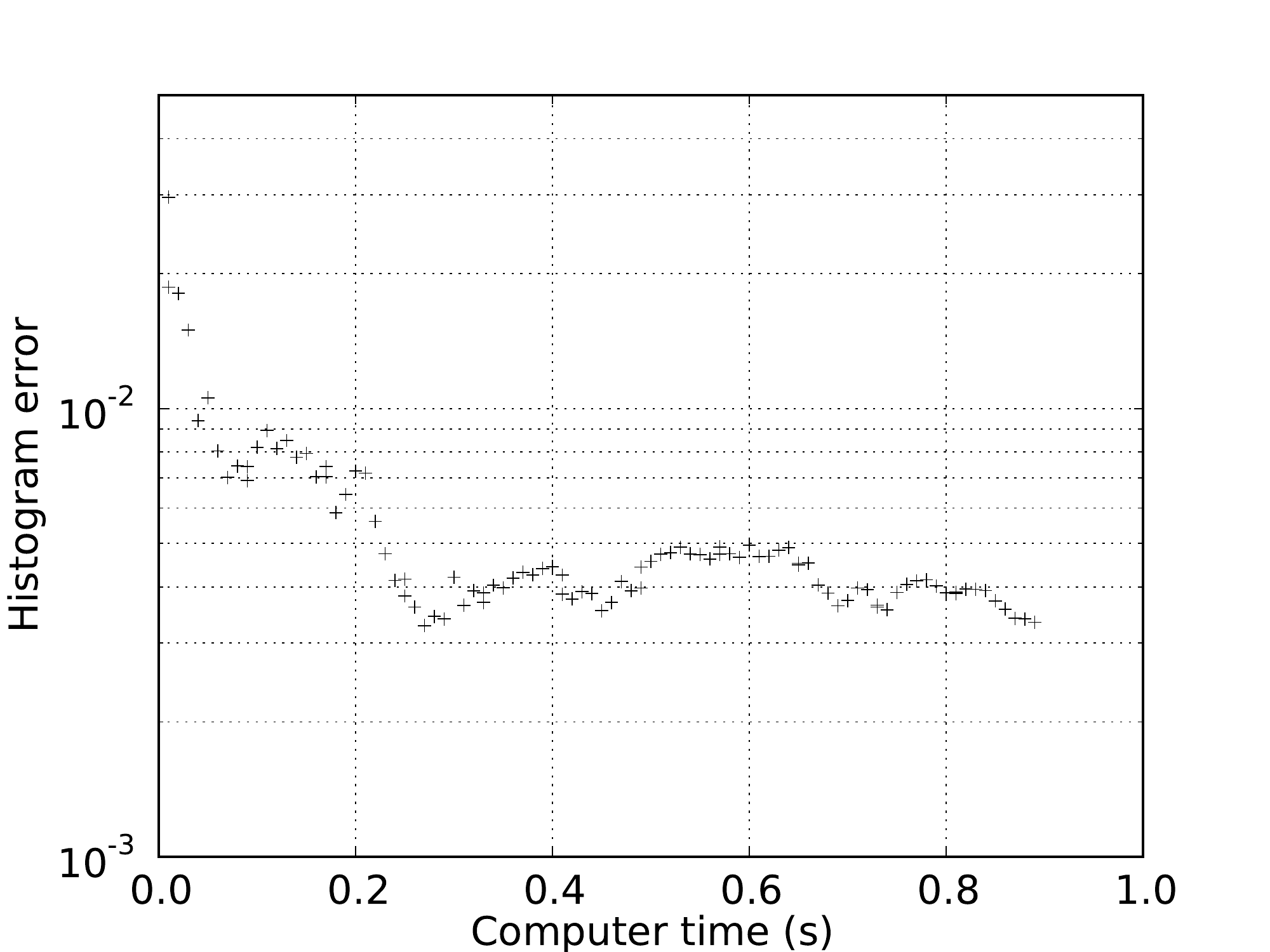}

}
\par\end{centering}

\begin{centering}
\subfloat[Local SF]{\includegraphics[width=6cm]{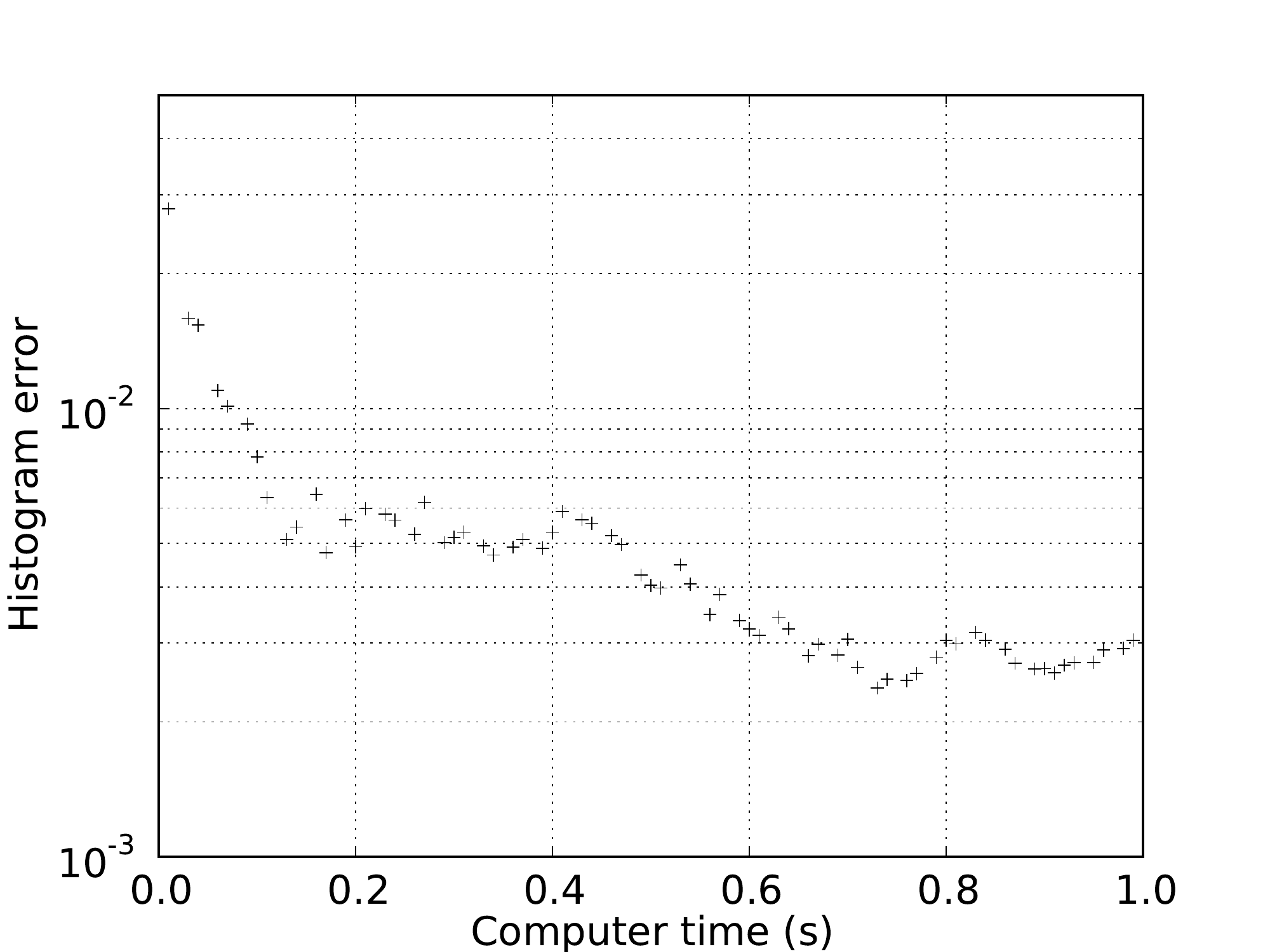}

}\subfloat[AISF]{\includegraphics[width=6cm]{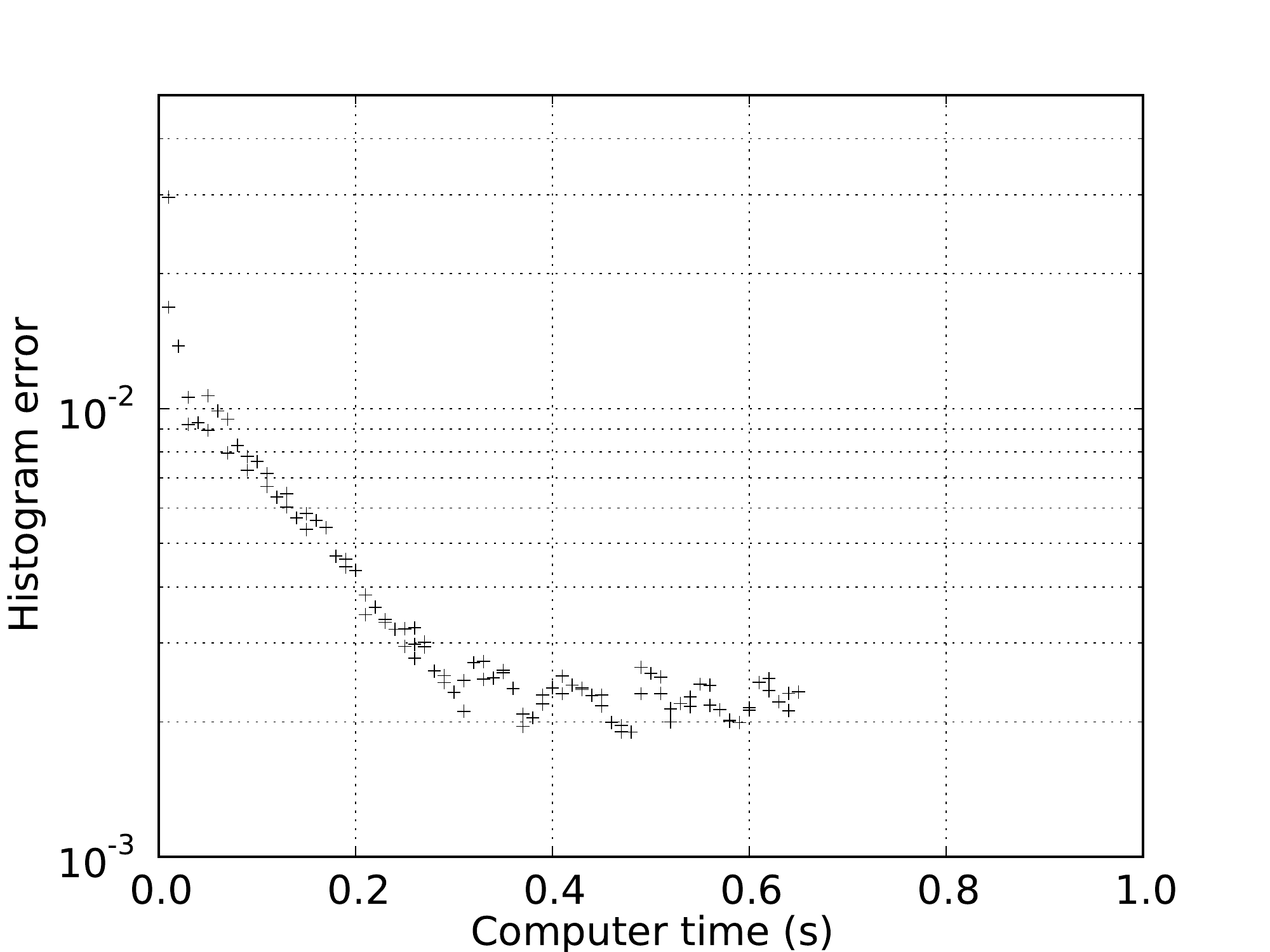}

}
\par\end{centering}

\caption{Histogram error versus computer time, for simulation from normal distribution.\label{fig:GaussConv}}
\end{figure}

A histogram of the simulated values is gathered, and compared with
a histogram calculated from the Gaussian distribution. We define the
histogram error using an $L^{\infty}$ measure, i.e. as the absolute
value of the maximum deviation from the exact result. The plots show
the histogram error as a function of simulation running time. The
simulations were run on one core of an Intel T4400 laptop processor,
but only the relative efficiencies matter here. Results are shown
in Figure \ref{fig:GaussConv} and are discussed below.

\subsubsection{Local MH}

Since a low correlation will be important for our later use on Levy
processes, we selected the transition probability width $0.75$, giving
the minimal correlation of approximatively 0.55, according to Figure
\ref{fig:CorrVsWidthLMH}. The variate generation rate was approximately
$3.4\cdot10^{6}/s$ in this case.

\subsubsection{AIMH}

As expected, as great improvement of the correlation was obtained
compared to local MH. Since the algorithm is more complicated, the
variate generation rate decreased to $2.3\cdot10^{6}/s$, or roughly
$60\%$ of local MH. However, the correlation was reduced to around
0.03. This dramatic decrease results in a faster histogram convergence
in terms of computer time. The lower correlation nature of this algorithm
will be noticeably beneficial when simulating Levy processes.

\subsubsection{Local SF}

We adjusted the basic SF algorithm parameters to give zero correlation,
and used a uniform transition probability density on $[-5,5]$. Thus
we are running it quite inefficiently as regards variate generation
speed, which turned out to be around $1.4\cdot10^{6}/s$. Despite
this, the histogram converges faster than local MH, due to the lack
of correlation.

\subsubsection{AISF}

As expected, when including the domain subdivision and adaptive algorithm,
the SF algorithm efficiency increased (in fact doubled) with a variate
generation rate of $3.1\cdot10^{6}/s$. The amount of correlation
here is still zero, which leads to the fastest histogram convergence.
So this is a great improvement at no cost other than increased algorithm
complexity. It is perfectly suited for simulating consecutive jumps
in jump-diffusion process paths.

\section{Jump-diffusion simulation\label{sec:CompSimJumpDiff}}

As previously mentioned, existence of correlations among jumps within
the same jump-diffusion path is disastrous. Let us check this by using
an local MH algorithm to simulate an exactly solvable stochastic process
\cite{Merton:1976} and compare distributions at $t=1$. It is defined
as
\[
X_{t}:=\sigma W_{t}+J_{t},
\]
\[
J_{t}=\sum_{i=1}^{N_{t}}Y_{i},\quad Y_{i}\sim N(\mu,\delta^{2}),\quad\mu,\delta\in\mathbb{R},
\]
where $N_{t}$ is a Poisson process of intensity $\lambda$. We choose
$\sigma=1$, $\mu=0$, $\delta=1$, and $\lambda=10$ in order to
get an appreciable number of jumps within the time domain $[0,1]$.

The exact PDF of this process for $t>0$ is \cite[Eq.(4.12)]{ContTankov:2003}
\[
p_{t}(x)=\exp(-\lambda t)\sum_{k=0}^{\infty}\frac{(\lambda t)^{k}\exp(-\frac{x^{2}}{2(t+k)})}{k!\sqrt{2\pi(t+k)}}.
\]

\begin{figure}
\begin{centering}
\includegraphics[width=8cm]{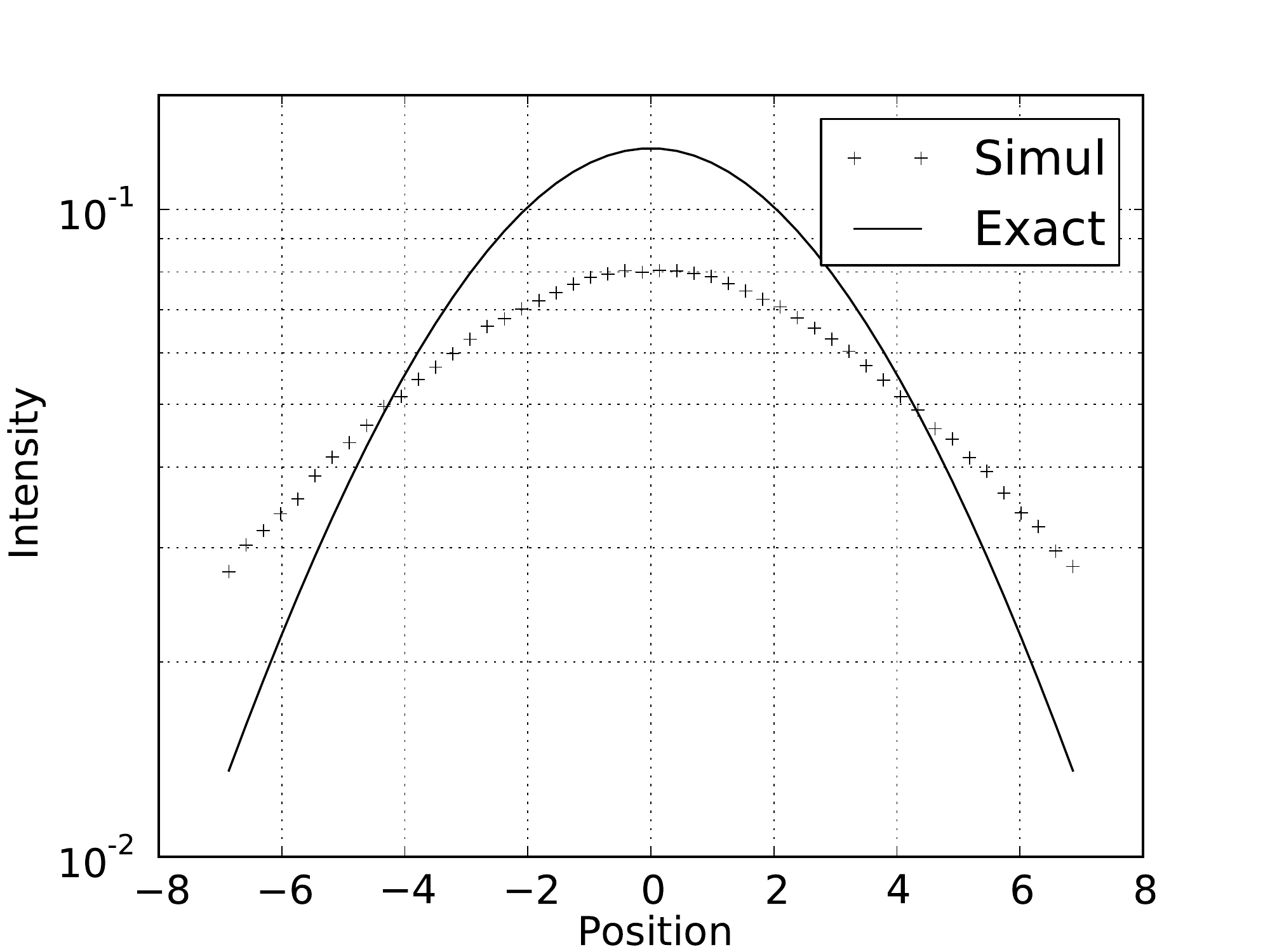}\caption{Local MH simulation of distribution at $t=1$ of Gaussian process,
and exact result. Due to jump correlations, the simulated histogram
has heavy tails.\label{fig:GaussProcLM}}

\par\end{centering}

\end{figure}

For local MH simulations, the positive correlations between subsequent
jumps leads to heavy distribution tails. Results confirming this are
shown in Figure \ref{fig:GaussProcLM}. The simulation of the jumps
used a localized uniform transition probability of width 4, which
resulted in a MH acceptance rate of around 0.56.

\begin{figure}
\centering{}\subfloat[AIMH]{\includegraphics[width=6cm]{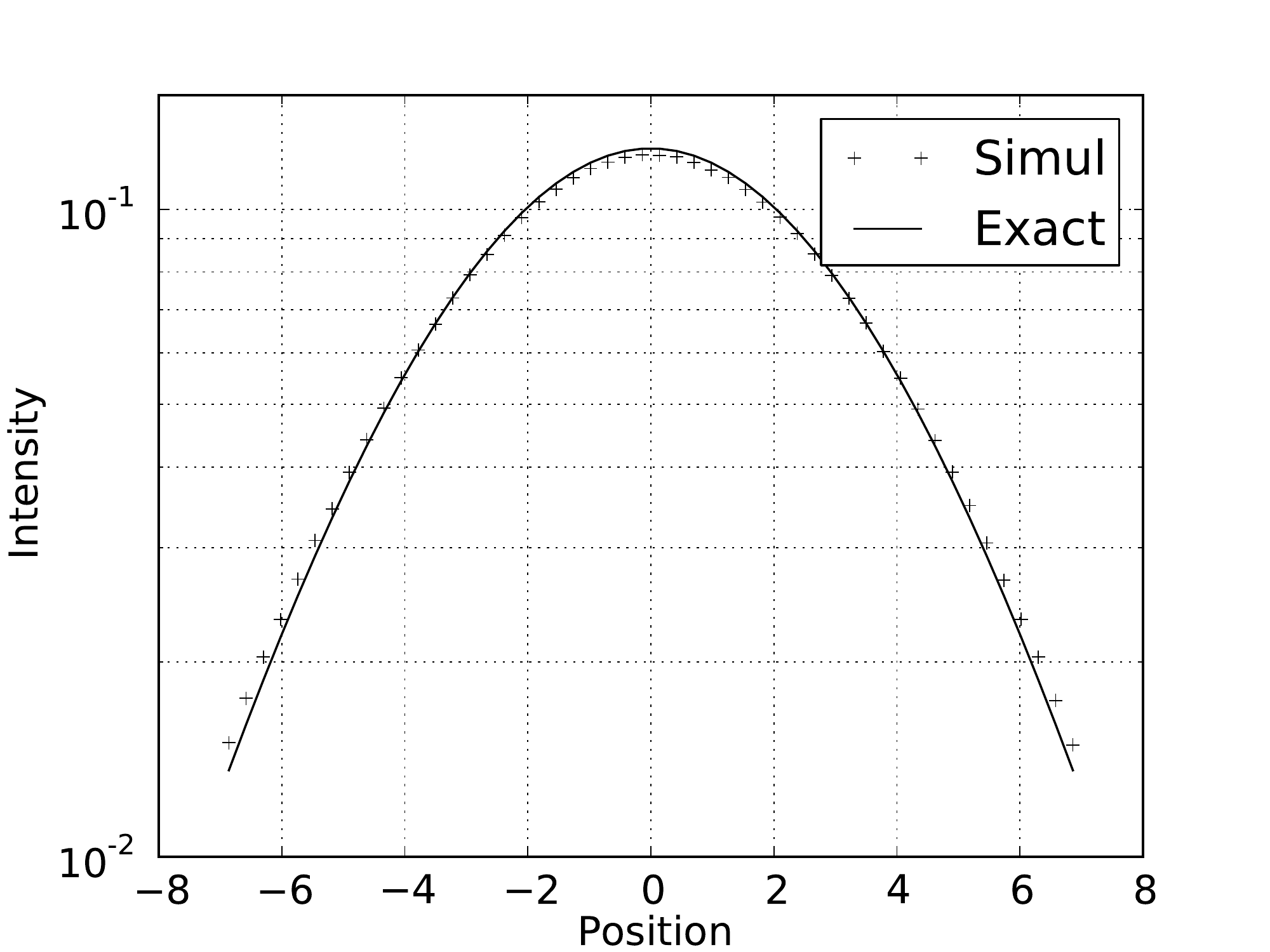}

}\subfloat[AISF]{\includegraphics[width=6cm]{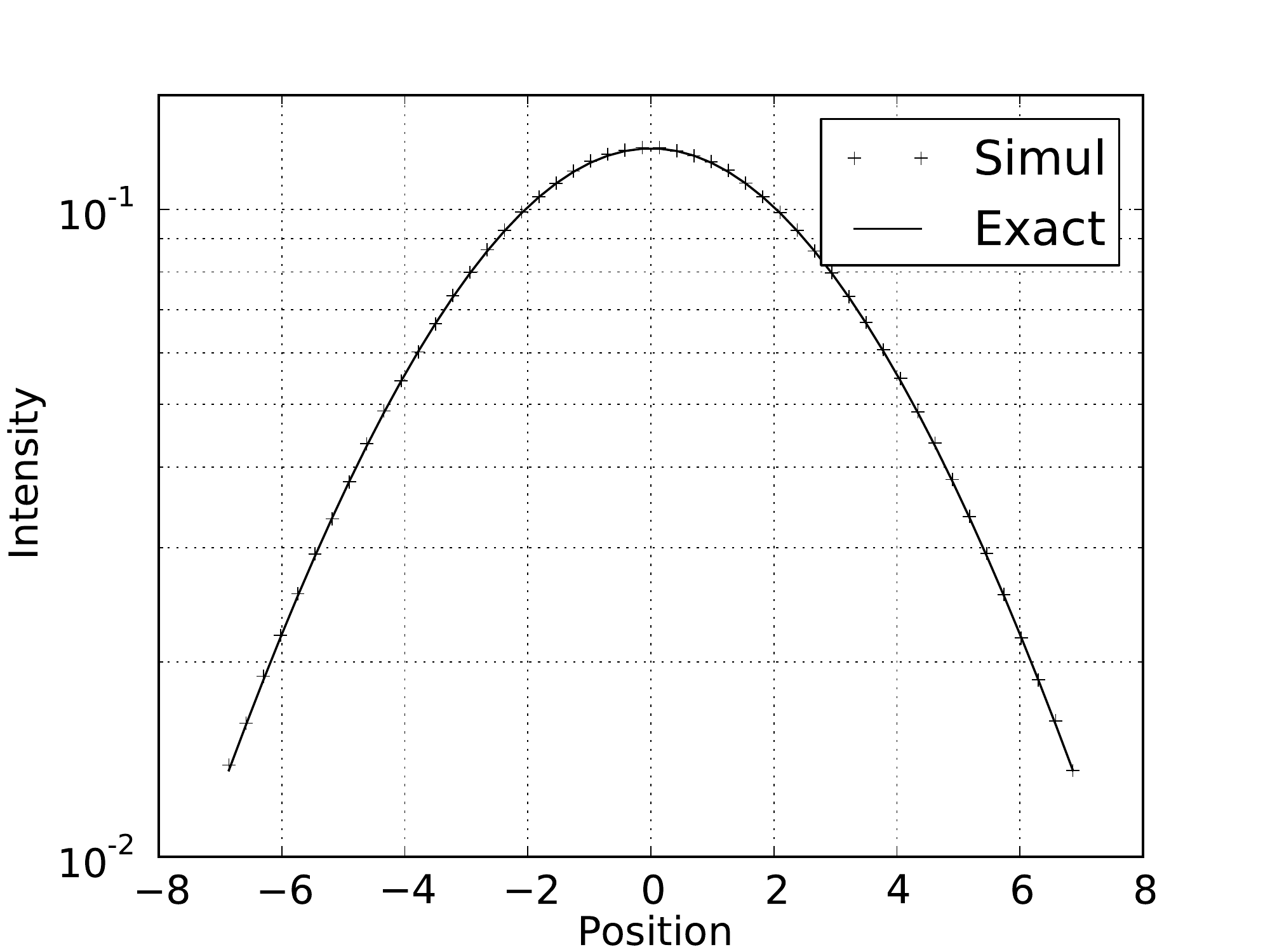}

}\caption{Distributions of Gaussian process at $t=1$ for AIMH and AISF simulations.
The continuous curve represents the exact result. Notice the slight
heavy tails in the AIMH case, due to nonzero jump correlations.\label{fig:DistGaussian}}
\end{figure}

We now turn to our serious attemps at jump-diffusion simulation. We
have simulated the same process using AIMH and AISF algorithms. In
this case we use the already known value of $\lambda=10$ in the simulations,
so the adaptability consists solely of the calculation of the discrete
subdivision distribution, and in the AISF case also of the local subdivision
supremum calculations. As opposed local MH, in this case the distribution
at $t=1$ approaches the exact curve, as seen in Figure \ref{fig:DistGaussian}.
Note that despite our use of upper/lower cutoffs in the implementation
of the discrete subdivision, the tails of the local MH  simulation
are still somewhat heavy.

\begin{figure}
\begin{centering}
\includegraphics[width=8cm]{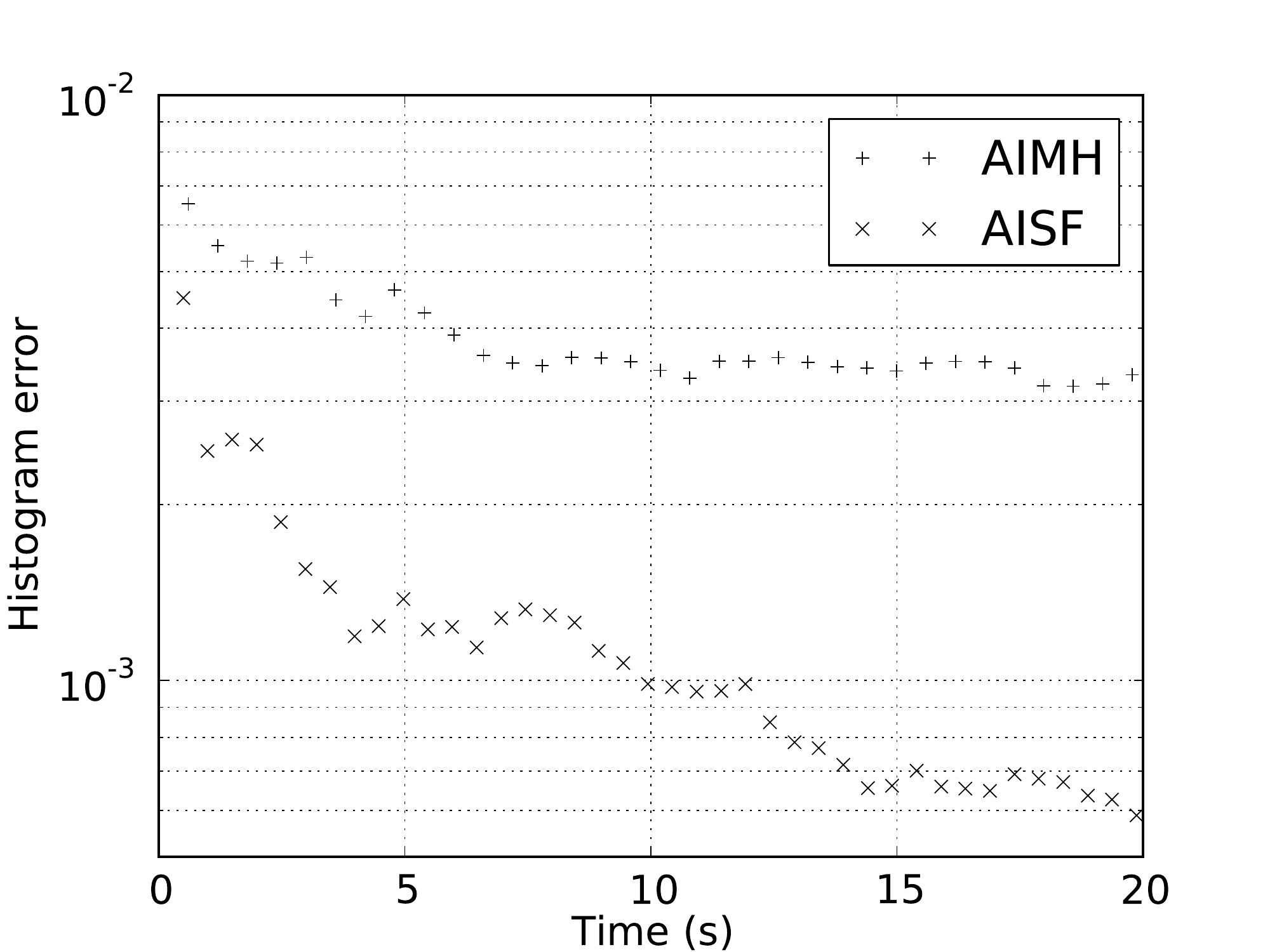}\caption{Error in distribution at $t=1$ for the Gaussian process, versus computer
simulation time. AISF gives the fastest convergence. The AIMH heavy
tails cause a lower bound on the error.\label{fig:ConvDistGaussian}}

\par\end{centering}

\end{figure}

We use the same $L^{\infty}$ measure of convergence detailed above.
Convergence results are shown in Figure \ref{fig:ConvDistGaussian}.
The AISF outperforms the AIMH algorithm due to its lower correlation
which leads to faster and more accurate convergence.

\section{Application to NIG and CGMY infinite activity pure jump processes\label{sec:SimInfAct}}

We present two applications within the realm of infinite activity
pure jump Levy processes, namely NIG and CGMY. To this end, we employ
the jump difffusion approximation of these processes \cite{AsmussenRosinski:2001}.

\subsection{Jump-diffusion approximation}

For an infinite activity Levy process, the intensity $\lambda$ is
undefined, since its defining integral \eqref{eq:intensity} diverges.
However, by virtue of \eqref{eq:LevyMeasCond}, it is possible to
approximate the infinitude of small jumps by a Brownian motion process.
Consider an infinite activity pure jump Levy process $L_{t}$. For
$\epsilon>0$, define the following subsets of $\Omega$,
\begin{eqnarray*}
A^{\epsilon,-} & := & \{|x|<\epsilon:x\in\mathbb{R}\}\\
A^{\epsilon,+} & := & \{|x|>\epsilon:x\in\mathbb{R}\}.
\end{eqnarray*}
These represent small and large jumps, respectively.

We now define the unique Levy measures $\nu^{\epsilon,-}$ and $\nu^{\epsilon,+}$
on these subdomains as follows. For any $\nu$-measurable $U\subset\mathbb{R}-\{0\}$,
define
\begin{eqnarray*}
\nu^{\epsilon,-}(U) & := & \nu(U\cap A^{\epsilon,-})\\
\nu^{\epsilon,+}(U) & := & \nu(U\cap A^{\epsilon,+}).
\end{eqnarray*}
Each of these Levy measures in turn define its own Levy process $L_{-}^{\epsilon}$
and $L_{+}^{\epsilon}$ of infinite and finite activity, respectively.
We have the unique process decomposition
\begin{equation}
L_{t}=L_{t}^{\epsilon,-}+L_{t}^{\epsilon,+}.\label{eq:InfActLevySplit}
\end{equation}
It follows from \eqref{eq:LevyMeasCond} that the intensity of the
large jump component process,
\[
\lambda^{\epsilon,+}:=\nu^{\epsilon,+}(\mathbb{R})
\]
is well-defined, and so $L_{t}^{\epsilon,+}$ is a compound Poisson
process. For small $\epsilon$, the small jump process $L_{t}^{\epsilon,-}$
can be often be approximated by the following Brownian motion \cite{AsmussenRosinski:2001}
\begin{equation}
L_{t}^{\epsilon,-}\approx\sigma(\epsilon)W_{t},\quad\sigma(\epsilon)^{2}:=\mathrm{Var}[L_{1}^{\epsilon,-}],\label{eq:BrownVol}
\end{equation}
where $W_{t}$ is a Wiener process. The variance must exist for the
approximation to be well-defined. A sufficient condition is \cite{AsmussenRosinski:2001}
\begin{equation}
\lim_{\epsilon\rightarrow0}\frac{\sigma(\epsilon)}{\epsilon}=\infty.\label{eq:AsmRosCondition}
\end{equation}
Thus we have the approximation
\begin{equation}
L_{t}\approx\sigma(\epsilon)W_{t}+L_{t}^{\epsilon,+}.\label{eq:BrownSmallJumpApprox}
\end{equation}
Any additional nonvanishing deterministic drift component of $L_{t}$
would appear trivially on the right hand side of this equation.

In these cases, the intensity $\lambda$ is not known beforehand,
and also depends on our choice of $\epsilon$. Since our algorithms
already calculate $\tilde{\nu}$, and the individual volumes of each
subset in the subdivision of $\Omega$ is known, it is easy to use
this to update an estimate for $\lambda$ concurrently.

\subsection{Simulation of NIG}

\begin{figure}
\begin{centering}
\subfloat[Cont/Tankov]{\includegraphics[width=6cm]{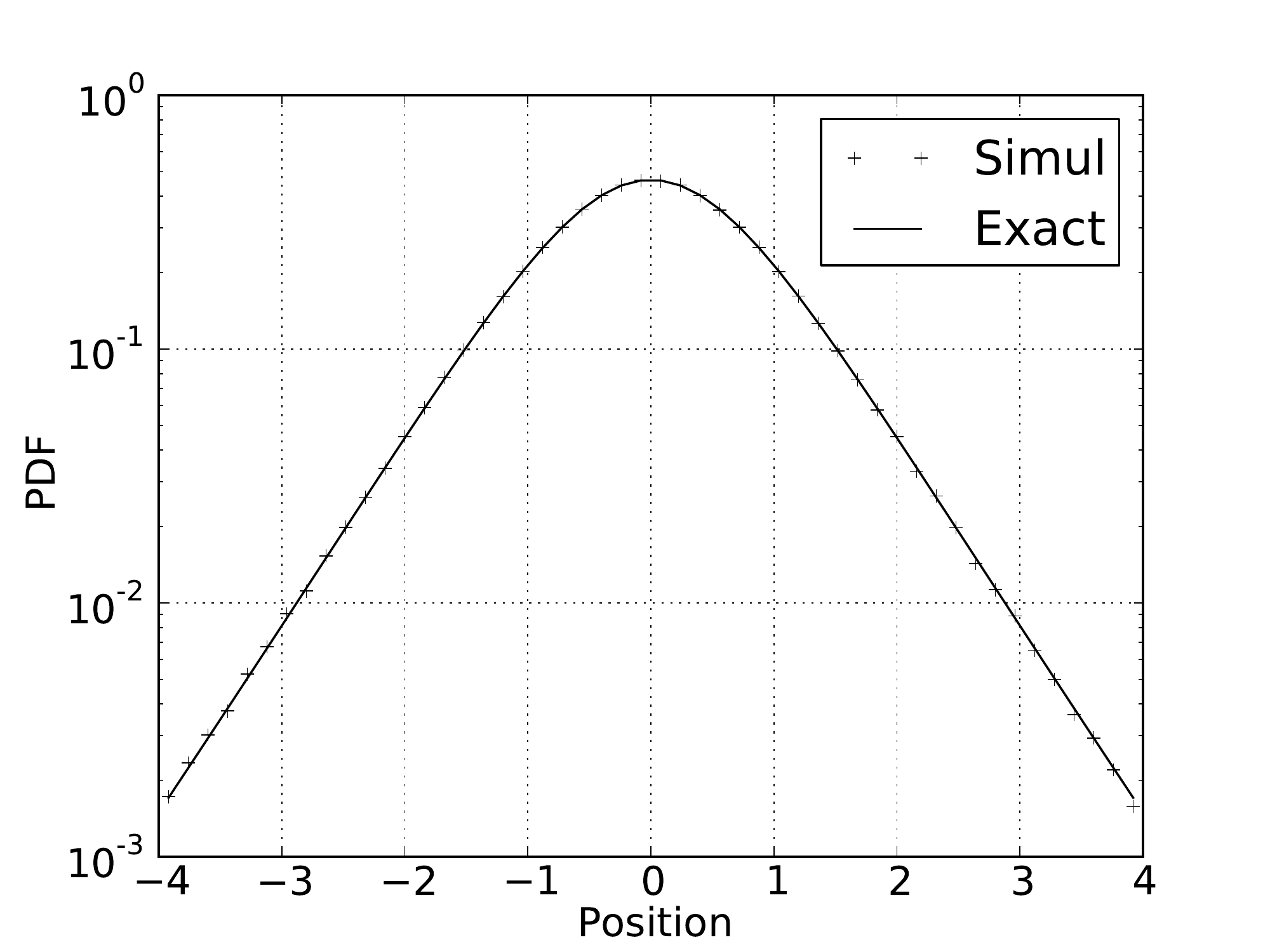}

}
\par\end{centering}

\begin{centering}
\subfloat[AIMH]{\includegraphics[width=6cm]{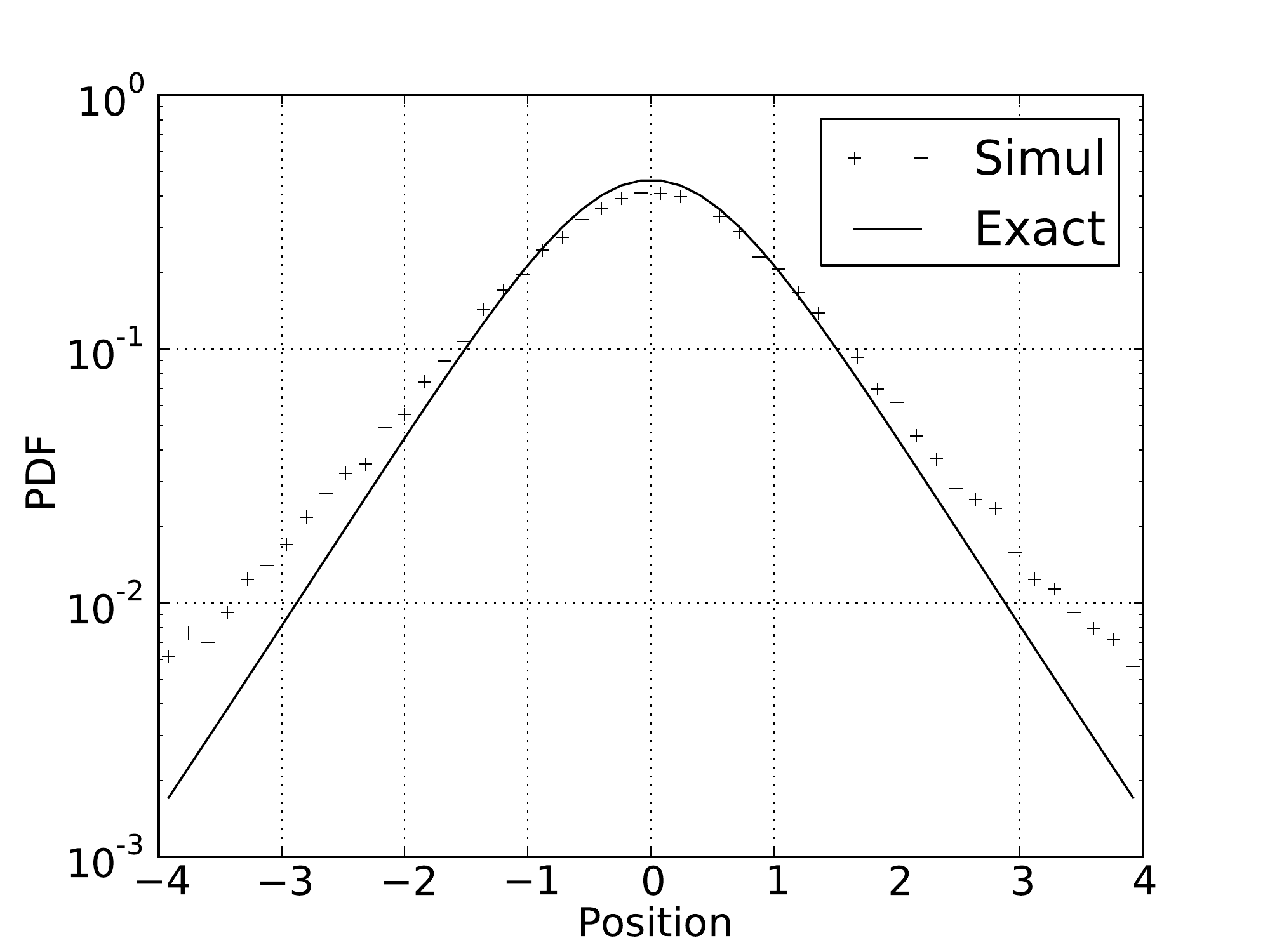}

}\subfloat[AISF]{\includegraphics[width=6cm]{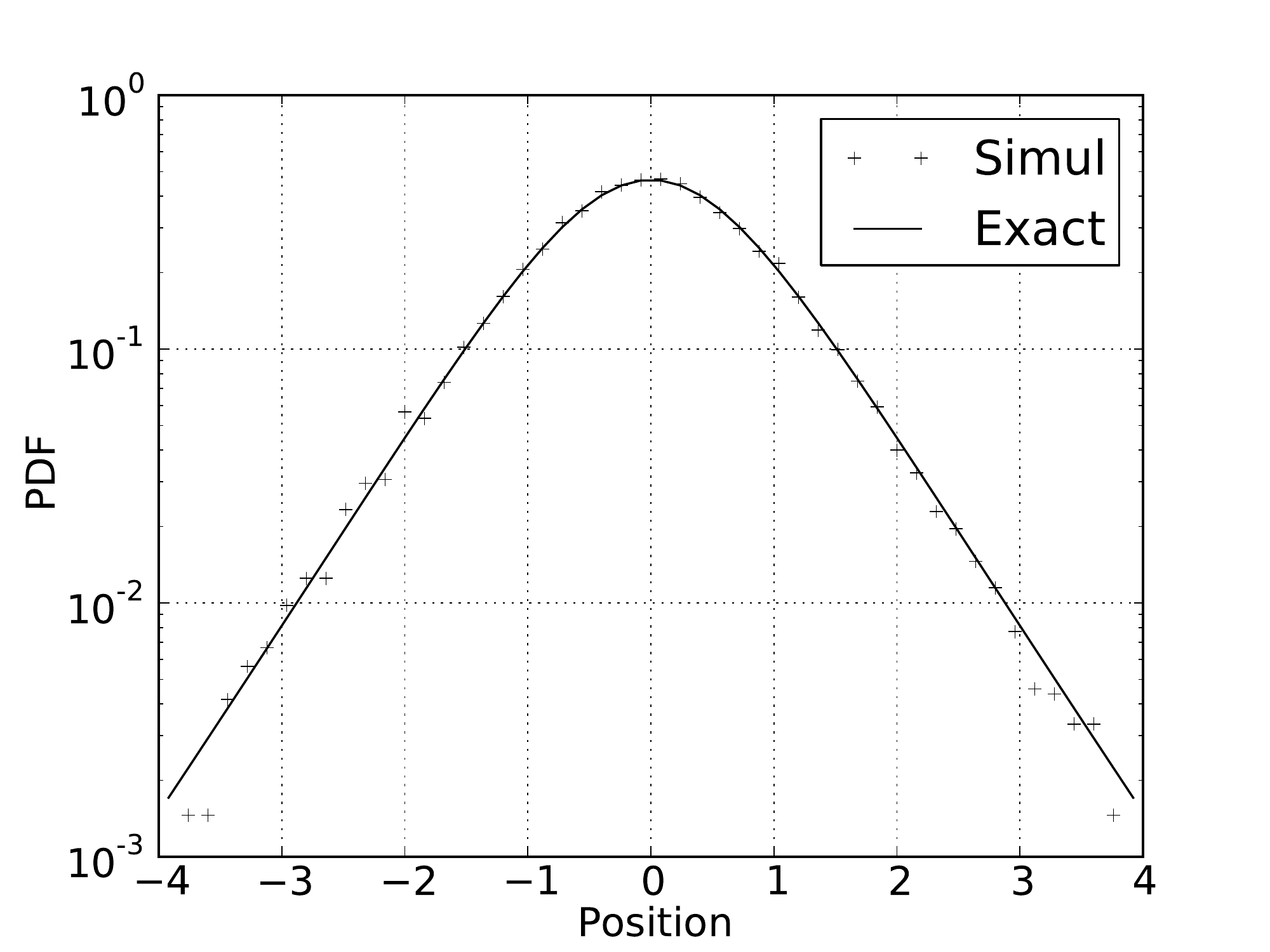}

}
\par\end{centering}

\caption{Simulated and exact NIG PDFs at $t=1$. The AIMH algorithm has heavy
tails.\label{fig:NIGSims}}

\end{figure}

We can now check the quality of the jump-diffusion approximation in
conjunction with our AIMH and AISF algorithms. Since the density of
NIG is known, and a direct algorithm \cite[Alg.6.12]{ContTankov:2003}
exists for simulating from its distribution at any time, we have ample
material to compare our results. The numerical results for the density
at $t=1$ are collected using logarithmic plots in Figure \ref{fig:NIGSims},
in order to emphasise the distributional tail behaviour. In terms
of the parametrisation used in \cite{ContTankov:2003}, we use $\sigma=1$,
$\theta=0$, $\kappa=1/2$. Figure \ref{fig:NIGSims} contains the
results from the direct simulation algorithm, as well as the results
of the jump-diffusion approximation using AIMH and AISF, where we
have used the small jump cutoff value $\epsilon=0.005$.

The volatility $\sigma$ of the Brownian component in the jump-diffusion
approximation, defined by \eqref{eq:BrownVol}, is calculated analytically
using an approximate expression for the Levy measure at small $|x|$.
The Levy jump density is
\begin{equation}
\nu(x)=\frac{\alpha\delta}{\pi|x|}e^{\beta x}K_{1}(\alpha|x|),\label{eq:NIGLevyDens}
\end{equation}
where $K_{1}$ is the irregular modified cylindrical Bessel function
of first order. As an asymptotic approximation of \eqref{eq:NIGLevyDens}
at small $|x|$, we use
\[
e^{\beta x}\approx1,\quad K_{1}(\alpha|x|)\approx\frac{1}{\alpha|x|}.
\]
Inserting the chosen parameter values, and using \eqref{eq:BrownVol}
which gives $\sigma^{2}$ as the second moment of the Levy jump density
on $[-\epsilon,\epsilon]$, we get
\[
\sigma\approx\sqrt{2\epsilon\delta/\pi}\approx0.067.
\]
Notice that this expression satisfies the Asmussen/Rosinski condition
\eqref{eq:AsmRosCondition}. Owing to this small value, the Brownian
part of the Levy process has a negligible influence on the distribution
tails at $t=1$, compared to the contributions from larger jumps.
\begin{figure}
\begin{centering}
\includegraphics[width=8cm]{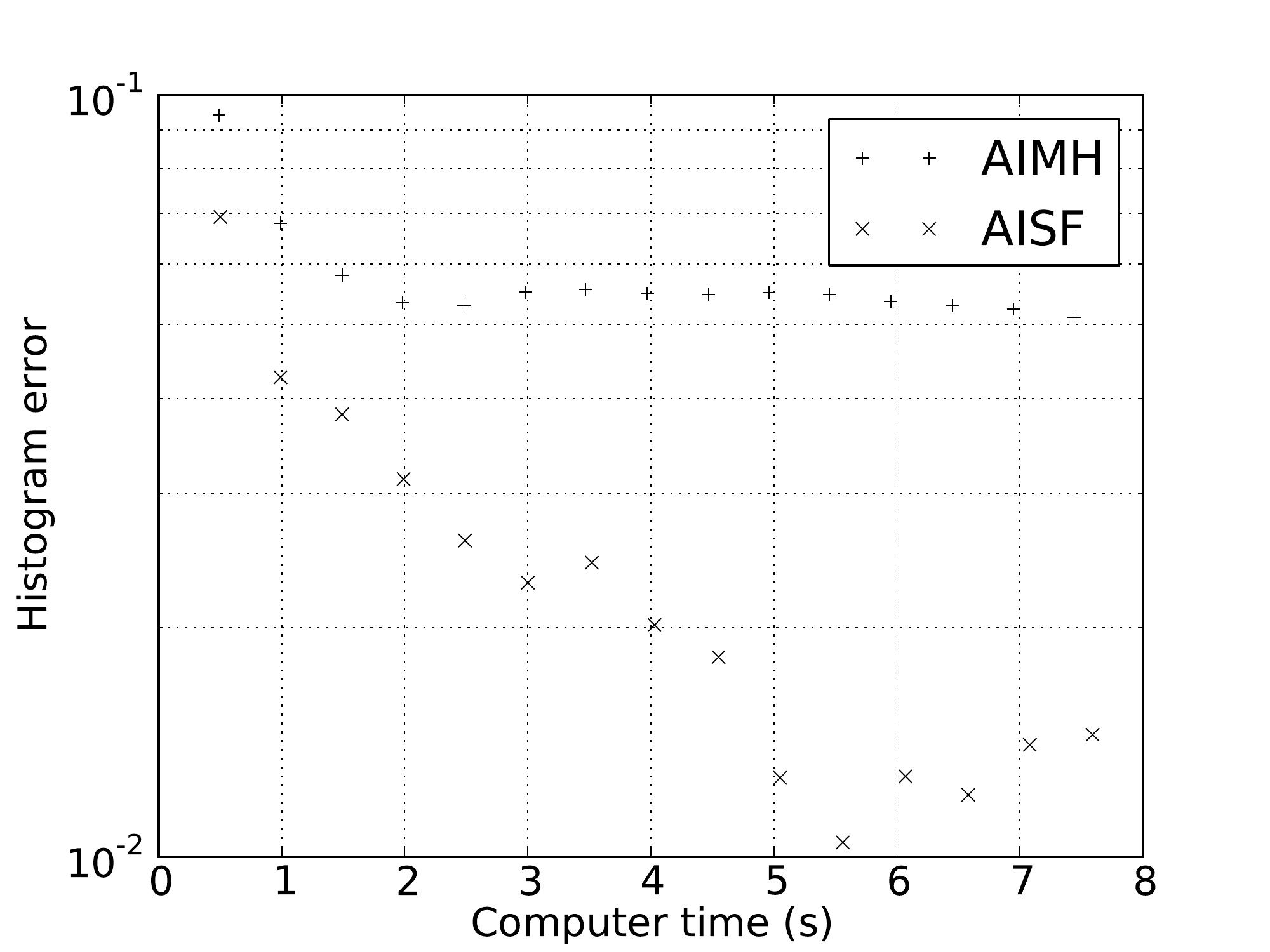}\caption{NIG histogram error at $t=1$, versus computer time. The AISF clearly
converges faster. AIMH reaches a steady state in the convergence plot
due to the error caused by the nonvanishing jump correlations.\label{fig:ConvNIG}}

\par\end{centering}

\end{figure}

The direct algorithm from \cite{ContTankov:2003} works best, as expected.
We are doing this to compare AIMH and AISF, for a known case with
an exact solution. The AIMH and AISF algorithms are general and can
be used on any Levy process that allows a jump-diffusion approximation.
This case gives us an indication of how trustworthy these methods
are when applied to more exotic Levy processes for which we do not
have an exact result or a simplified simulation methods as in this
case.

Most sources of errors are common to both algorithms. These are related
to inaccuracy in the calculation of the $\sigma$ and $\lambda$,
the latter being related to the cutoff imposed on the jump domain
$\Omega$. This cutoff will cause an inaccuracy in $\lambda$ since
the total weight of $\Omega$ will not be calculated. Improved subdivision
schemes of $\Omega$ are possible, e.g. employing coordinate transformations
that transform $\Omega$ into a compact domain. We have not done this,
since it unimportant for the matters we are considering.

We note that AIMH will never be completely correlation-free, and will
therefore tend to produce heavy tails as is obvious in Figure \ref{fig:NIGSims}.
No such problem exists for AISF. The cutoff on $\Omega$ will naturally
lead to weak tails, as seen in the plot. This can be remedied by enlarging
the cutoff value, and/or treating large jumps differently.

The measurement data for convergence analysis is plotted in Figure
\ref{fig:ConvNIG}, from which one readily sees that the AISF algorithm
converges faster and more exactly.

\subsection{Simulation of CGMY}

As a second example of a pure jump infinite activity process, we turn
to CGMY \cite{Carr01stochasticvolatility}. In this case we have no
closed form expression for the distribution. We do however have the
characteristic function, from which the distribution can be obtained
by means of a numerical inverse Fourier transform. We have performed
this calculation, and used the result as a the benchmark against which
our path simulation algorithms are compared.

In this case, the Levy density is
\begin{equation}
\nu(x)=\left\{ \begin{array}{ll}
Ce^{-Mx}/x^{1+Y} & ,x>0\\
Ce^{-G|x|}/|x|^{1+Y} & ,x<0.
\end{array}\right.\label{eq:CGNY_LevyDens}
\end{equation}
 The parameter space for the CGMY model is $C,G,M>0$ and $Y<2$.
By expanding the Levy density in negative powers of $x$ around the
origin and keeping only the lowest order terms, we get
\[
\sigma(\epsilon)\sim\epsilon^{1-Y/2},
\]
so by \eqref{eq:AsmRosCondition}, the jump-diffusion approximation
is valid only for $0<Y\leq1$. In fact, the volatility can in this
case be expressed exactly in terms of the incomplete gamma function,
which we will not show here.

We simulated the CGMY process with parameters $C=G=M=1$, $Y=1/2$,
using $\epsilon=0.005$, and 100000 paths. In this case, the volatility
for our choice of parameter values is
\[
\sigma\approx0.022.
\]
Also in this case it has negligible influence on the distribution
tails at $t=1$.

\begin{figure}
\begin{centering}
\subfloat[AIMH]{\includegraphics[width=6cm]{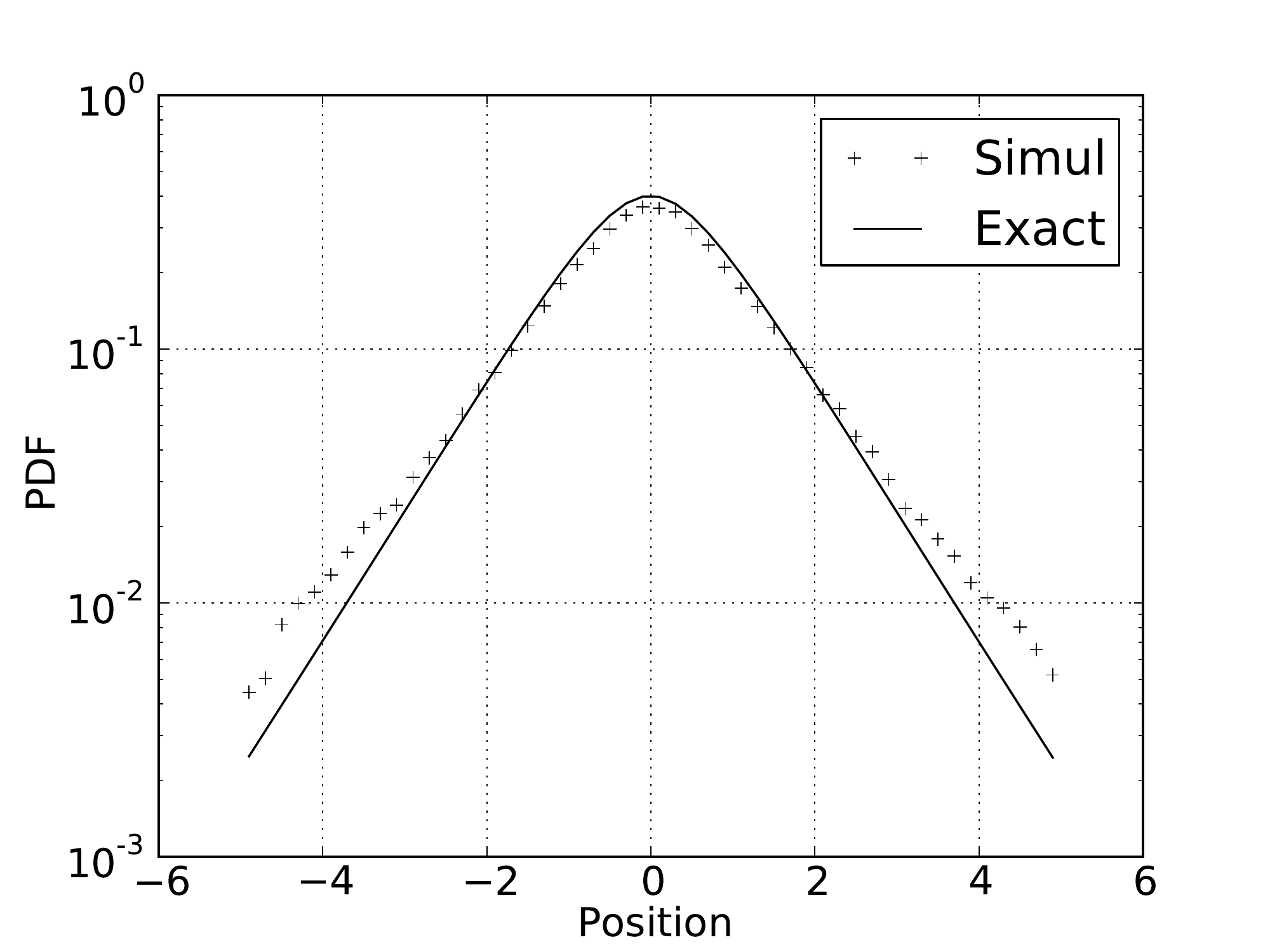}

}\subfloat[AISF]{\includegraphics[width=6cm]{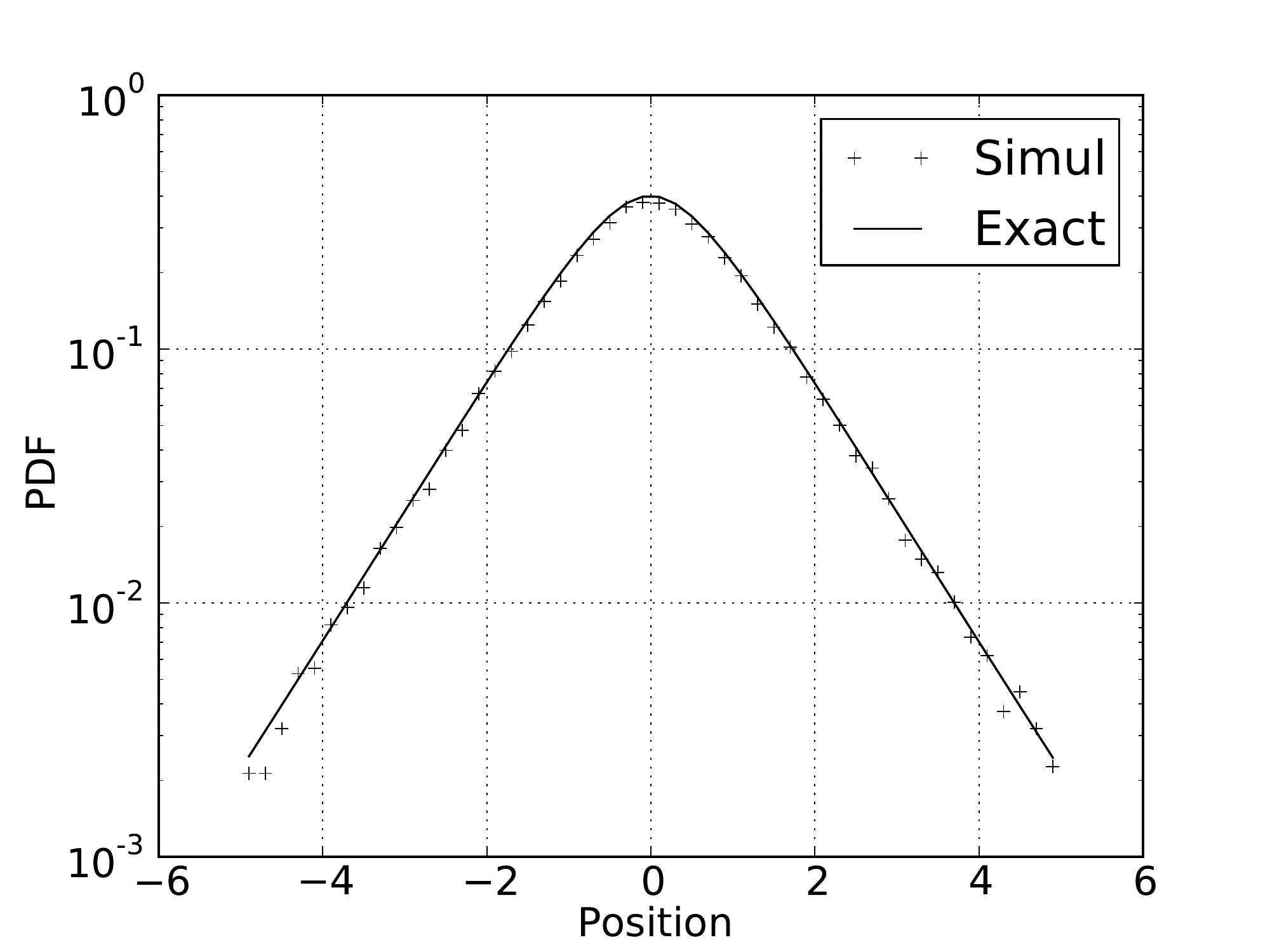}

}
\par\end{centering}

\caption{CGMY PDF at $t=1$, simulated and exact results. Note again the heavy
tails in the AIMH case.\label{fig:CGMYDensity}}
\end{figure}

\begin{figure}
\centering{}\includegraphics[width=8cm]{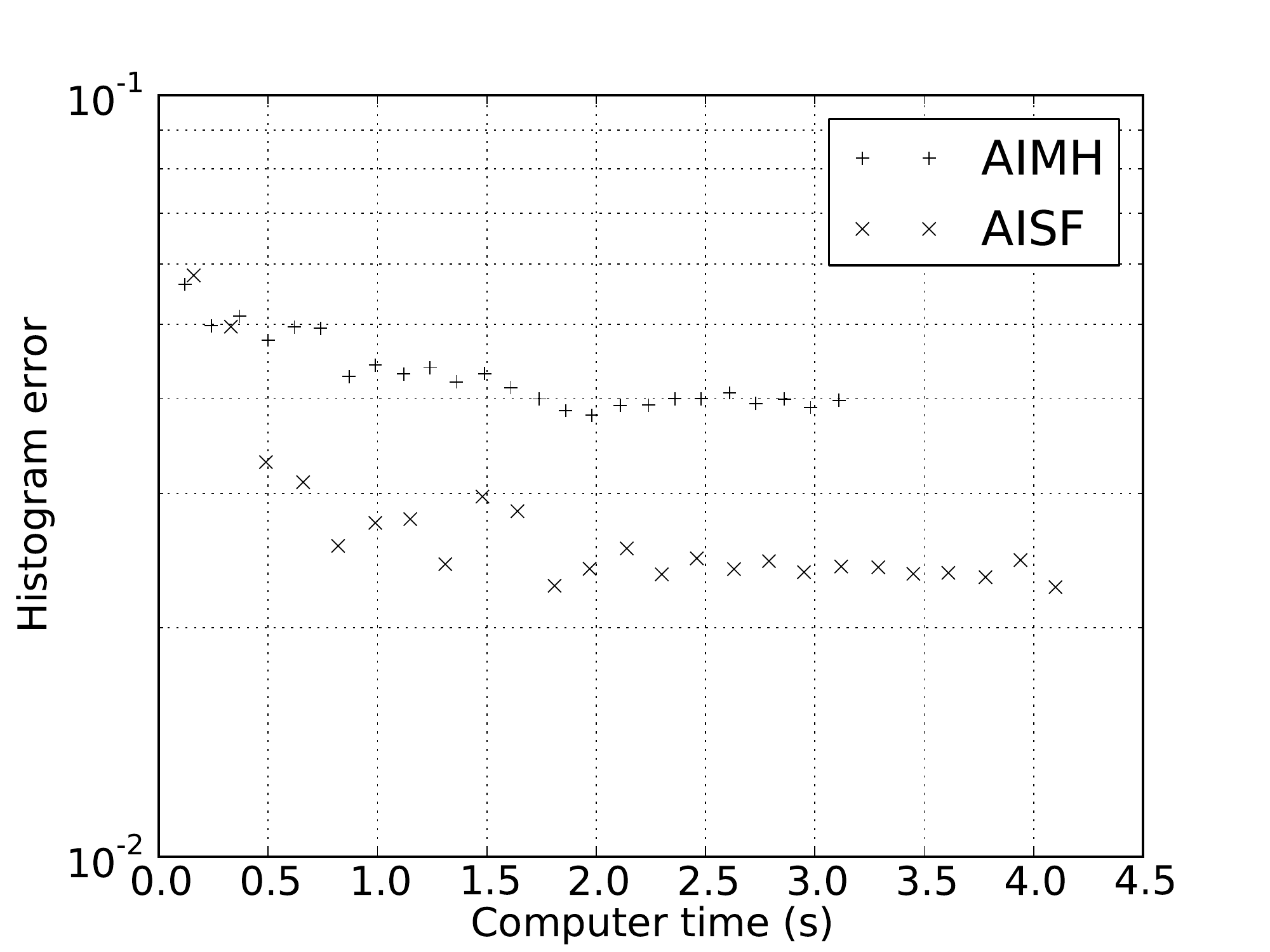}\caption{CGMY histogram error at $t=1$, versus computer time.\label{fig:ConvCGMY}}
\end{figure}

The results for the distribution at $t=1$ of the process is given
in Figure \ref{fig:CGMYDensity}, and the convergence measurements
are plotted in Figure \ref{fig:ConvCGMY}. The conclusions are similar
to the NIG case.

\section{Conclusions}

We have studied two algorithms for jump-diffusion and infinite activity
pure jump process simulation via jump-diffusion approximations. Most
importantly, we have studied the proposed SF algorithm based on a
stochastic step function. It has some advantages over MH accept/reject
algorithms. It is possible to configure it to have an arbitrarily
small autocorrelation. Our simulations show that in the case of simulation
of Levy processes, this algorithm can represent an improvement over
the MH method that we have considered. The numerical results show
an improvement in the tails of the distribution of the Levy process
at a given time while at the same time converging faster. The MH algorithm
will tend to give heavy tails, due to the problem of positive correlations
between large jump values.

In our computer simulations, we also implemented a subdomain discretization
with a corresponding adaptively improved discrete probability distribution.
This method helps to reduce correlations for the MH algorithms, since
the subdomains themselves are drawn without using an accept/reject
algorithm. In the SF case, it improves the variate generation speed
while maintaining the correlation-free nature of the method.

\bibliographystyle{siam}
\bibliography{references}

\end{document}